\documentclass[12pt]{amsart}
\usepackage{graphicx} 
\title{Accessibility for the Working Mathematician}
\author{Julius Ross}
\usepackage{hyperref}

\newcommand{\TeXml}{\TeX ML}

\usepackage{amsthm,amsmath,amssymb} 
\usepackage{tikz}
\newtheorem*{WCAG}{WCAG Success Criterion}

\theoremstyle{definition}
\newtheorem{exercise}{Exercise}[section]

\newtheorem*{digression*}{Digression}
\newtheorem*{examplequestion}{Example Exam Question}

\begin{document}

\maketitle


\emph{If you only take one thing away from this document I would like it to be this:  Creating accessible documents requires authors have a working understanding of accessibility, and write with accessibility in mind.}

\section{Introduction}

This is version 2 of this document, written in April 2026.

\subsection{Disclaimer}

I claim no deep expertise on the subject of accessibility, and there will likely be errors and misinterpretation in what I have to say.  In particular none of this should be taken as legal advice, nor as assurance that if you follow this document you will meet any specific accessibility standards or satisfy any particular legal regulations.  At the time of writing version 2 of this document, I serve as Head of the Department of Mathematics, Statistics, and Computer Science at the University of Illinois Chicago and as Chair of the American Mathematical Society (AMS) Advisory Group on Accessibility.   However, this document is written purely in my capacity as a Professor of Mathematics, errors are my own and the opinions expressed do not represent the AMS, my department or my university.

\subsection{About}

This document is an informal attempt to share knowledge I have gathered about writing accessible mathematical documents that I hope may be of use to others.   It is aimed at mathematicians or scientists taken in the broadest sense, whether they be teachers or researchers, who are interested in writing accessible electronic documents that contain mathematical formulae.  I do not try to be comprehensive, and although parts focus on electronic documents produced in a higher education setting in the United States, the majority will be relevant to those working in other countries or settings.

Given the above disclaimer,  it may be helpful to think about whether my writing such a document has any value at all.  I argue that it does for the following reasons:

\begin{itemize}
\item Existing accessibility documentation is often geared at general audiences, and mathematical writing has its own conventions and needs (not least our need for formulae, some of which may be long or complicated).
\item There are various authoring software options available and in development, each with advantages and disadvantages.  Making definitive recommendations is thus very difficult, opinions differ, and the best choice likely depends on context.
\item Accessibility documents and standards can be long, technical and are unlikely to make it to the New York Times best-seller list.  I attempt to focus on the parts that are likely to be of most relevance to mathematical documents.
\item Various accessibility best practices and compliance requirements are subjective.  Just as the mathematical community has, over time, developed standard conventions and notations for mathematical writing we may find we want to develop conventions and notations for accessible mathematical writing.  Since this process is ongoing, an informal discussion is more appropriate than a formal list of recommendations.
\end{itemize}

\subsection{Acknowledgments}

Much of what is included in this text I have learned directly from others.  I thank Loretta Bartolini and David Jones from the American Mathematical Society and Peter Krautzberger not only for sharing their comments and expertise, but also for the production and public release of the {\TeXml} tools. I also thank Deyan Ginev for answering my questions about the {\LaTeX ML} project.   

I thank my colleagues on my campus' Digital Accessibility Steering Committee as well as Nick Christo, Michael Gintz, Sophia Hamilton and JaEun Jemma Ku for discussions and feedback. I also thank my colleague Fred Drueck for discussions on the topic of accessibility and for technical support.

Finally I thank my colleagues on the AMS Accessibility Advisory Group, Jorge Balbas, David Farmer, Nicola Poser, Jan Reimann, Daniel Thompson as well as Tyler Kloefkorn for their efforts in this arena; many of the ideas in the revisions that appear in the second version of this document come from them.

\section{What is Accessibility?}

\subsection{Accessibility}

Accessibility refers to the proactive steps that ensure that a wide range of individuals, including those with disabilities, have independent access without the need for modifications.   Accessibility can apply in many situations, such as to building access, to employment processes and, as is the focus here, to electronic documents.

It is worth emphasizing that creating accessible electronic documents has advantages that go beyond serving people with disabilities.  For instance, the addition of captions on videos can help those who need to work in noisy environments.  Also, accessible features in electronic documents may allow them to be read on a wider range of devices, such as phones, which could be important for those on low income who do not own a device with a larger screen. Furthermore, some accessibility features facilitate machine-reading of documents, which can have advantages for indexing and searching.

\subsection{Accessibility $\neq$ Accommodations}\label{sec:accommodations}

Accessibility and accommodations are related concepts, but they serve different purposes in ensuring inclusivity for people with disabilities.  

As said above, accessibility (sometimes referred to as ``universal access") refers to the proactive steps that ensure that a wide range of individuals have independent access without the need for modifications.   Accommodations, on the other hand, are specific adjustments or modifications made to meet the individual needs of a person with a disability when an environment or service is not inherently accessible.

Those teaching in higher education will likely be familiar with the idea of accommodations in the classroom setting.  Although regulations may differ from country to country (or state to state), students have the right to receive ``reasonable accommodations" in the pursuit of their studies.  In practice this often means institutions have systems that give students the opportunity to disclose disability status, and the accommodations required to meet these needs are then made on a case-by-case basis by administrators or instructors of the courses in which the student is registered.  Similar accommodations systems are often available for employees.

Accessibility and accommodations require a slightly different mindset in that accessibility demands anticipation of those with different needs.  Meeting accessibility in educational documents reduces the burden on students with disabilities, gives access to those with undisclosed or unknown disabilities, and serves students who may be engaging with the material without formally registering.  For example, course descriptions, syllabi of courses, and sample examinations all need to be written in an accessible manner so that all students can make their course choices independently, which typically happens before the formal process of accommodations within a class can take place.

\subsection{Accessibility $\neq$ Compliance}\label{sec:compliance}

Because the notion of accessibility is open-ended, it is often convenient, or necessary, to have specifications or standards that precisely describe certain accessibility features.  For example, for accessibility to physical spaces these standards are often included as parts of building codes.

For electronic documents there are a number of different accessibility standards.    It is useful to note from the outset that even full compliance with any of these does not necessarily equate to full accessibility.   It is not hard to find examples of documents that are fully compliant but still can improve their accessibility and, vice versa, examples of documents that may fail to meet a particular specification but are still accessible to a wide range of people.  

There may be situations in which an author is most interested in accessibility and may choose to use such standards as a (likely very good) guide to best-practice.  On the other hand, there may be different situations in which an author is more interested in doing precisely what is required to remain compliant with a certain standard. 

\subsection{Web Content Accessibility Guidelines (WCAG)}
Among all the accessibility specifications for electronic documents the only one we consider here is the Web Content Accessibility Guidelines (WCAG) \cite{WCAG}.   There are different versions of these (e.g. WCAG2.0 and WCAG2.1) and within each version there are different ``levels" of compliance, level A, level AA, level AAA, which are progressively stricter.  

When discussing some specific aspects of accessibility below I will cite some parts of these guidelines, but the reader wanting to be fully compliant will need to reference the full guidelines, which contain not only the specifications but also explanations and criteria for success and failure.  Sometimes it will be convenient to paraphrase or quote just a snippet of the WCAG which will be indicated with the $\simeq$ or $\subset$ symbol.  Also I emphasize that we will be looking at just a selection of the WCAG that are likely to be most relevant to electronic mathematics documents; the full guidelines have much more with pieces that apply to user-interface, audio, video, etc.

\subsection{US Legal Framework}\label{sec:DOJ}
Disability rights are covered in many legal frameworks.  In the United States, higher education institutions have obligations under Title II of the Americans with Disabilities Act (ADA), that, in particular, gives students with disabilities the right to accommodations in pursuit of their education.  In 2024, the United States Department of Justice gave a final ruling under Title II of the ADA requiring state and local government websites and mobile apps to conform to WCAG specifications.

We will not delve into the details of this ruling here, but to say that for state institutions serving a population of 50,000 persons or more it requires all electronic documents, whether they be aimed at students, staff or visitors, be compliant to WCAG 2.1 level AA by April 2027 (with a later date specified for smaller state institutions).   It is worth pointing out also that some states have their own regulations regarding accessibility.

A summary of the requirement is publicly available \cite{dojsummary} as is the full ruling \cite{dojfull}.  It applies broadly not just to websites but also to electronic documents, and although there are some exemptions (for instance for archived content) these exemptions are narrowly defined.  In particular this ruling applies to all internal and external websites, homework systems as well as electronically distributed lecture notes, examinations, syllabi etc.    I believe it also applies to research documents.

Although such rules are established at the federal or state level, they are implemented and interpreted at the institutional level, meaning different institutions may apply them somewhat differently.   It is up to mathematicians and others to work with their institution's compliance office (or disability resource center as appropriate) to explain their needs and help shape their institutions interpretation of these rules.  In my opinion, the best approach prioritizes accessibility over compliance, acknowledges that context matters, and accepts that we can only do what current technology allows.  

This compliance date was initially set for April 2026, and the extension to the new date of April 2027 was given with just a few days notice.   In my view this extension offers an opportunity for individuals and units to develop their accessibility frameworks, learn about accessibility, understand what tools are available and give constructive feedback to their institution's compliance departments as to the needs of mathematicians.  Hopefully there will continue to be technological advances that will facilitate writing STEM rich accessible documents in this period.

It remains to be seen what the future of this new ruling will be, including how it will be interpreted, or how it will be enforced and on what timeline. For many institutions, compliance is a significant challenge, and while most larger institutions appear to have made substantial and meaningful progress, smaller ones may lack the campus resources to navigate this ruling in the same way.  

\section{Key concepts}\label{sec:keyconcepts}

\subsection{Writers First, Technology Second}

Creating accessible documents requires authors to have a working understanding of accessibility, and write with accessibility in mind\footnote{Readers who are keen to make highly accessible documents with minimal author intervention or accessibility knowledge are recommended to look into PreTeXt \cite{PreTexT}.}

Said another way, authors should consider that their audience may include people with various needs, who require documents written to be compatible with the accessibility software tools they use. (Unfortunately, for many mathematics documents it is not clear the author has considered their audience at all, but that is a different story!)

As a trivial example, on virtually any software it is possible to compose red text on a green background.  This is likely a bad idea for various reasons, not least as it makes such text invisible to the approximately 4\% of the population who are color-blind.  Other examples may be less obvious and are the subject of the following topics that we will discuss in subsequent sections:

\begin{itemize}
\item Alternative Text
\item Alternative Text for Mathematical Formulae
\item Alternative Text in Examinations
\item Navigation and Headings
\item Links
\item Colors
\item Fonts
\end{itemize}
Meeting the requirements within these topics to make documents accessible can often be achieved rather easily, simply by changing a few writing habits.  

It may be in the near or distant future that there are AI tools that will do a lot of this for us, but my use of such tools have given me mixed results.  I expect this is because so often accessibility decisions need to be made in the context of the text and intended audience.  Furthermore current AI systems can suffer from ``hallucinations", so it is important that the user has a good idea of what good accessibility looks like.

There exist various tools that claim to test if a document or website is accessible and/or compliant, but these tools seem to simply flag if a document contains material that fails compliance in some of the most common ways (e.g.\ choice of colors).  Such checkers are certainly a good place to start, and can be excellent triage tools, which can be especially useful when looking at a large number of documents.  I have not seen a tool that would give me confidence that a document is truly accessible or fully compliant, especially when it comes to mathematical documents.





\subsection{Alternative Text}\label{sec:alternativetext}

Accessible documents need to be able to be read by screen readers, and so steps must be taken with regard to non-text content such as images, graphs, diagrams, etc.  The WCAG says:

\begin{WCAG}[$\simeq$ 1.1.1 Non-Text Content]
All non-text content that is presented to the user has a text alternative that serves the equivalent purpose [with some specific exceptions].
\end{WCAG}

\begin{figure}
\begin{center}
\includegraphics[alt={A torus in three dimensions as an example of an image that could have multiple alternative texts},width=10cm]{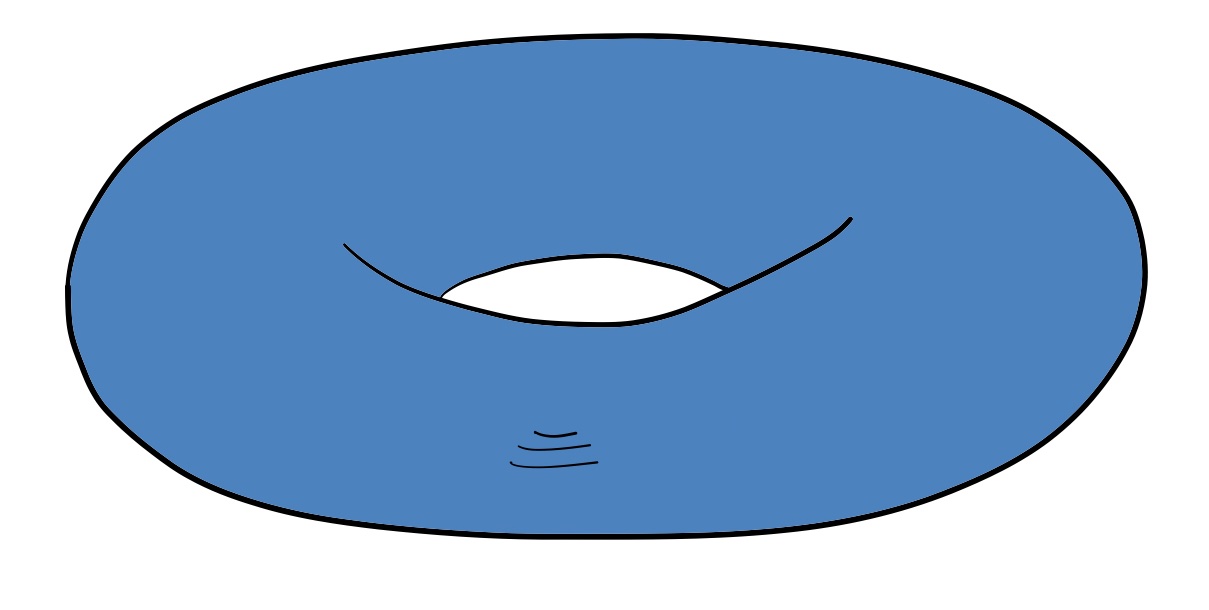}\end{center}
\caption{A torus in three dimensional space}\label{fig:torus}
\end{figure}

Notice here the requirement that this text serve ``equivalent purpose", which means that the same image or graph could have different alternative texts in different contexts.   Consider for instance the image in Figure \ref{fig:torus}.  This could be trying to illustrate different concepts, for instance the graph of an implicit function, a torus sitting in three-dimensional Euclidean space,  a surface of revolution, the gluing of a cylinder, a flat torus, a general genus one Riemann surface etc.\   Some of these may be mathematically equivalent, but a useful alternative text has to depend on the context and the intended audience.

There are various best-practices for writing good alternative text, and further resources can be found at \cite{daigramcenter}.    The following are adapted from \cite{alttext,daigramcenter}:

\begin{itemize}
\item Alternative text should be short and communicate the same information as the visual content. 
\item Alternative text should describe the content of the image, not describe how it looks.
\item Alternative text should not contain unnecessary information, or repeat information already in the text or caption.  In particular the alternative text should not start with ``An image of" or similar unnecessary text. 
\item For complex images (for instance graphs with data) alternative text should be short and descriptive and point to the location of a longer text that is comprehensive \cite{WIAcomplex}.

\end{itemize}

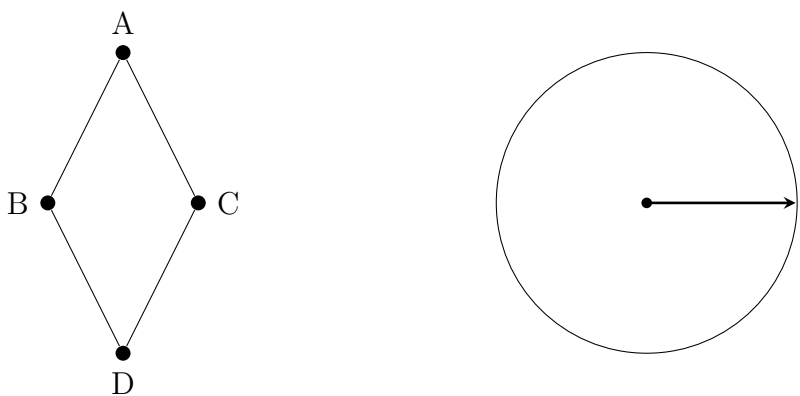
\begin{figure}
\begin{minipage}{0.45\textwidth}
\begin{center}
\begin{tikzpicture}
    \node[circle, fill=black, inner sep=2pt, label=above:A] (top) at (0,2) {};
    \node[circle, fill=black, inner sep=2pt, label=left:B] (a) at (-1,0) {};
    \node[circle, fill=black, inner sep=2pt, label=right:C] (b) at (1,0) {};
    \node[circle, fill=black, inner sep=2pt, label=below:D] (bottom) at (0,-2) {};

    \draw (bottom) -- (a);
    \draw (bottom) -- (b);
    \draw (a) -- (top);
    \draw (b) -- (top);
    \end{tikzpicture}
    \end{center}
\end{minipage}
\hfill \begin{minipage}{0.45\textwidth}\begin{center}
\begin{tikzpicture}
    \draw (0,0) circle (2);
    
    \fill (0,0) circle (2pt);
    
    \draw[->, line width=1pt, >=stealth, scale=1.5] (0,0) -- (1.32,0);
\end{tikzpicture}
\end{center}
\end{minipage}\caption{Figures that can appear in mathematical documents}\label{fig:exercise}
\end{figure}

\begin{exercise}
For each image in Figure \ref{fig:exercise} write down a few different contexts in which this image could appear in a mathematical document, and for each such context write a suitable alternative text.
\end{exercise}

\subsection{Alternative Text for Mathematical Formulae}\label{sec:alttextformulae}

To be accessible, formulae and equations need to be readable by screen-readers, and the part of the WCAG concerning alternative text applies since formulae are considered ``non-text".   As we will see below, it is not recommended that authors (only) use alternative text for formulae since there is widely used technology that allows for much better accessibility.

For simple equations such as $$f(x)=2$$ a suitable alternative text would be ``\emph{f of x equals two"} or perhaps ``\emph{f open parentheses x close parentheses equals two}".  For an expression such as $$V^{\perp}$$ perhaps a good alternative text is ``\emph{V perp}" or ``\emph{The orthogonal complement of $V$}" if that is what the superscript actually means in this case.  From this latter example we already see that the most useful alternative text depends on context, and is place where authors would have to use their judgment if they were writing them manually.

Consider now a more complicated formula, such as the following:
\begin{equation}X = (\zeta-\lambda)n + \eta + \frac{\lambda - \zeta - \mu -\nu}{2} + \frac{\eta n - \frac{1}{4}(\lambda - \zeta +\mu-\nu)^2-\theta}{X'}\label{eq:euler}\end{equation}

Any linear alternative text to this is going to be long, and thus much less comprehensible.  For imagine trying to understand this equation purely from the following text:\medskip

\emph{``X equals open parentheses zeta minus lambda close parentheses times n plus eta plus a fraction whose numerator is lambda minus zeta minus mu minus nu and whose denominator is 2 plus a fraction whose numerator is eta times n minus a quarter open parentheses lambda minus zeta plus mu minus nu close parentheses squared minus theta and whose denominators is X primed."}\medskip

Clearly such an alternative text is too long to be useful, so we need something better. Thankfully there are various tools that allow a vision impaired reader to explore such a formula with keyboard strokes, thereby focusing on pieces of such a formula (such as particular numerators or denominators) and have screen readers read and reread those pieces (I am referring here to MathJax and/or MathML, but will not go into the details of what these are).   Most of the commonly used authoring tools in STEM now have the ability to make equations accessible using these technologies without additional author intervention.

\begin{digression*}
The formula in \eqref{eq:euler} is taken from Euler's \emph{Opvscvla Analytica}, published in 1783 which was several years after he had turned completely blind.  Despite this, Euler apparently became even more productive and was able with the aid of scribes to produce approximately one mathematical paper per week.
\end{digression*}

\subsection{Alternative Text in Examinations}\label{sec:alttextexams}

Providing alternative text can pose a challenge for images used in tests or examinations since giving too much information in a question can render it useless (and remember alternative text is available to all readers, not just those with accessibility needs).  

\begin{figure}
\begin{center}
\begin{tikzpicture}
    \draw[->] (-0.5,0) -- (6.5,0) node[below] {\small $x$};
    \draw[->] (0,0) -- (0,4) node[left] {\small $f(x)$};
    \draw[very thin, gray] (0,0) grid (6,3);
    \draw[domain=1:5, smooth, samples=100, thick, blue] 
        plot (\x, {1/3*(\x)^3 - 3*(\x)^2 + 8*(\x)-4});
   \foreach \x in {1,2,3,4,5} {
        \node[below] at (\x,-0.2) {\small \x}; 
    }
\end{tikzpicture}
\end{center}
\caption{Graph Associated to Example Exam Question}\label{fig:exam}
\end{figure}
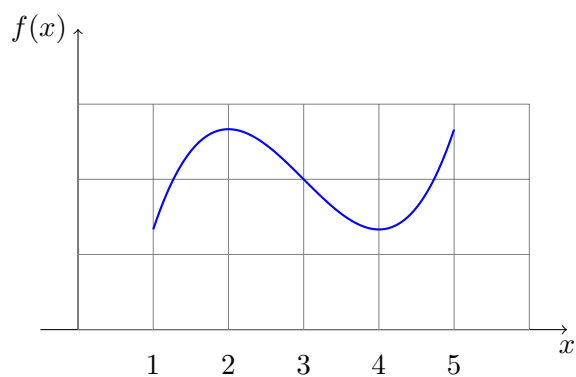

For example consider the following exam question: 

\begin{examplequestion}
\emph{Find all the local minima and local maxima of the function $f$ illustrated in Figure \ref{fig:exam} and explain your reasoning}. 
\end{examplequestion}

It is not clear to me that there is a good alternative text for this graph that does not invalidate the exam question.  The best I can come up with is something like ``\emph{The graph of a continuous function $f(x)$ from $[1,5]$ to $\mathbb R$ that strictly increases between $x=1$ and $x=2$, strictly decreases between $x=2$ and $x=4$ and strictly increases between $x=4$ and $x=5$}".  But this alternative text is not ideal as it gives clues that the correct answer involves the numbers $2$ and $4$, involves the concept of continuity that might not have been covered in the class, and uses phrases that are within the expected answer!\medskip

On this issue the WCAG say:

\begin{WCAG}[$\subset$ 1.1.1 Non-Text Content]
Test: If non-text content is a test or exercise that would be invalid if presented in text, then text alternatives at least provide descriptive identification of the non-text content.
\end{WCAG}

My interpretation (and opinion) is that authors have options here.  One option is to only use questions that are not invalidated by alternative text, which should be preferred when it does not diminish the purpose of the exam.

A second option is to do as the specifications say, and simply provide an alternative text that gives a ``descriptive identification".  In the example above this descriptive identification could be something like \emph{``The graph of a function whose purpose is to test students' knowledge of minima and maxima".}  One may argue that this second option might be compliant but is not accessible (c.f.\ Section \ref{sec:compliance}).  This could be true, but it is important to remember two things.  First, even this alternative text has value, in that somebody using it will understand the specific purpose of this figure as well as what skill or knowledge the exam question is trying to assess.  Second,  accessibility is not the same as accommodations (c.f.\ Section \ref{sec:accommodations}) so students who cannot use the graph as given should have accommodations to have their knowledge assessed another way.

\subsection{Navigation and Headings}\label{sec:navigation}

Accessible documents need to be structured with titles and headings in a way that makes them operable in many ways, such as only with a keyboard.    Relevant WCAG success criteria include:


\begin{WCAG}[$\subset$ 2.1.1 Keyboard]
All functionality of the content is operable through a keyboard interface [...]
\end{WCAG}
\begin{WCAG}[1.3.1 Info and Relationships]
Information, structure, and relationships conveyed through presentation can be programmatically determined or are available in text.
\end{WCAG}

\begin{WCAG}[2.4.2 Page Titled]
Web pages have titles that describe topic or purpose.
\end{WCAG}
\begin{WCAG}[2.4.6 Headings and Labels]
Headings and labels describe topic or purpose.
\end{WCAG}

Authoring tools nearly always have methods to include titles, headings and labels and the above success criteria will likely be met as long as these methods are used correctly.  For example a document should have a top level title that uses the title feature, followed by a heading at level 1 and using nested headings at lower levels as appropriate.   These titles and headings should be meaningful in that they help navigate around the document (a good practice even outside of accessibility considerations).  For instance in a document with many theorems, the heading ``Proof" is less meaningful than the heading ``Proof of Theorem A".  In fact accessibility best-practice is that each heading's text should be unique \cite{twingenuity}.

Skipping the title is not accessible, and skipping heading level 1 and just using heading level 2 is less accessible than not doing so (skipping the title is not WCAG compliant because titles can appear in various places outside of the document, such as search results or window title bars).  If authors want a title or heading to be at a particular font size they should use the correct title or heading level and adjust the document styling appropriately.   

 Authors should not inadvertently hide headings and labels using font styles.  For instance using a paragraph text but set in slightly larger font and/or bold may give the same visual appearance as a heading, but is not accessible as navigation to this heading will not be picked up by screen readers.

Finally authors should refrain from using other methods to organize content on the page, or indicate titles.  For instance centering text is purely visual, and unlikely to be picked up by screen readers.  And tables should only be used for displaying data, and never for organizing content on a page.
\subsection{Links}

\begin{WCAG}[$\subset$ 2.4.4 Link Purpose In Context]
The purpose of each link can be determined from the link text alone [...]
\end{WCAG}

To meet this criterion, authors should refrain from using the url in the link text, or from having link text separated from the description.  A compliant example is \href{https://pi.math.cornell.edu/~hatcher/AT/ATpage.html}{Algebraic Topology by Allen Hatcher}.   Two non-compliant examples are 

\noindent \href{https://pi.math.cornell.edu/~hatcher/AT/ATpage.html}{https://pi.math.cornell.edu/\textasciitilde hatcher/AT/ATpage.html} 
or text of the form ``Allen Hatcher's Topology book that can be found \href{https://pi.math.cornell.edu/~hatcher/AT/ATpage.html}{here}."

\subsection{Colors}

What is commonly referred to as ``red-green colorblindness" confuses different conditions and masks the fact that there are many ways that people perceive color.  (For those interested, there are ways to explore how different color palettes may be perceived by different groups \cite{coloringforcolorblindness}.)  

Since it is not possible to fully know the needs of your audience in this regard, the WCAG says the following:

\begin{WCAG}[1.4.1 Use of Color]
Color is not used as the only visual means of conveying information [...] or distinguishing a visual element.
\end{WCAG}

So authors should use information in addition to color, such as shape or text or a legend to convey meaning.  For example a figure with the graph of two functions should not only color alone to indicate which function is which, but instead label the graphs or indicate the difference in some other way.

Related to color is the topic of contrast, meaning that when multiple colors are used (such as a different color for text and background) a palette is chosen with colors that are sufficiently distinct.  

\begin{WCAG}[$\simeq$ 1.4.3 Contrast]
The visual presentation of text and images of text has a contrast ratio of at least 4.5:1 [with exceptions for large text, incidental text and logotypes].
\end{WCAG}

\begin{WCAG}[$\simeq$ 1.4.11 Non-text Contrast]
The visual presentation of [Graphical Objects] have a contrast ratio of at least 3:1 against adjacent colors.
\end{WCAG}

The technical meaning of this ratio can be found within the WCAG (and was chosen to give access to readers with 20/40 vision which is typical vision of elders around the age of 80).  Thankfully there are easy to use resources that provide color palette that meet this requirement (e.g. \cite{accessiblepallete}).  The criteria for graphical objects would apply, for instance, to colors used in histograms or pie charts.

\begin{figure}\begin{center}\includegraphics[alt={Two pictures of the same tiger in the undergrowth. 
 In picture (a) the tiger appear camouflaged as to it has colors similar to the undergrowth whereas in picture (b) the tiger is much more visible as it has orange-black colors that contrast with the undergrowth},width=13cm]{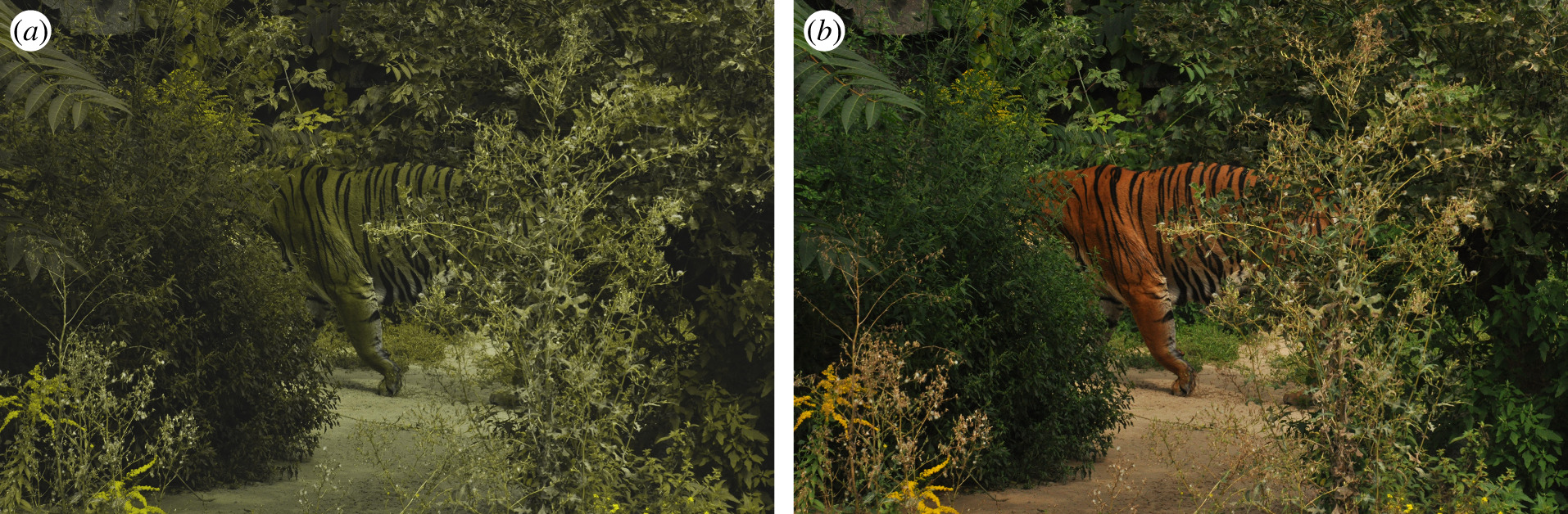}\end{center}
\caption{\small Image from \cite{tigers} of a tiger from the point of view of a simulated (a) dichromat and (b) trichromat. }\label{fig:tigers}

\end{figure}

\begin{digression*}
 A tiger's red-orange striped fur may not appear to most humans as an optimal form of camouflage.   But their typical prey (e.g.\ deer, boars) are dichromats which gives them a form of red-green color blindness, making the tigers appear to them as blending into the undergrowth (see Figure \ref{fig:tigers}).  An example from nature where it matters to think of your audience!
\end{digression*}

\subsection{Fonts}
It is well understood that certain fonts can be more or less legible and there are fonts specifically designed for the vision-impaired or those with dyslexia \cite{fontsforaccessibility}.
There are relatively few fonts available with good mathematics support and so, for better or worse, mathematicians (by which I really mean {\LaTeX} users) commonly choose between the same few serif fonts (e.g.\ Latin Modern or STIX).  It is worth pointing out that Google's Noto Sans is a sans-serif font was written with legibility in mind and has mathematics support.

Among accessibility best-practices is that electronic documents  allow the user to change style and/or font (e.g.\ the WCAG talks about the ability to change the size and spacing of a font).  There are good accessibility reasons to want to give the user the ability to change font completely which, if implemented, may also give authors more freedom in their font choice.   How easy this is to achieve will depend on the authoring software used; this is likely possible with HTML output but may be harder with some PDF output (see further discussion below).

\begin{digression*}
Generally speaking, writing accessibly means letting go of micromanaging certain aspects of your document.  One of these is in the difference between ``script" and "calligraphic" characters (which are often available in {\LaTeX} using the ``\textbackslash mathscr" and ``\textbackslash mathcal" commands).  My understanding is that current technology considers these as ``variants" of the same letter \emph{with the same meaning} \cite{variants}, and as such it is a challenge for the distinction to be accessible to those using screen readers.  

Authors are still free to choose one variant or the other, but should not using both variants for the same letter in a given document unless it is very clear from context what is meant.  So for instance if you end up proving ``\textbackslash mathscr A= \textbackslash mathcal A" then your choice of notation has gone wrong.
\end{digression*}

\begin{digression*}
Technology built to improve reading accessibility is not new.  Linn Boyd Benton and  Morris Fuller Benton were font designers circa 1900, who put effort into understanding the interplay between legible fonts and eyesight.  Among their creations was the Century Roman family of fonts made for legibility, a variant of which is now available with mathematics support under the name TeX Gyre Schola. 

I like to think that the Bentons would appreciate the recent technology advances on electronic accessibility.  Linn Boyd had failing eyesight later in life and apparently shortly before his death received a patent for ``for an important improvement in the larger printing types used in newspaper headings" \cite{bentons}.
\end{digression*}




\section{Choice of Authoring Tools}\label{sec:tools}  

Accessible electronic documentation requires at least three different kinds of tools that need to be compatible with each other: file formats, authoring software and reading/viewing software. 

\subsection{File Formats}

The following is a list of the main document formats I have seen used for mathematical documents, each of which offer various accessibility features:

\begin{itemize}
\item HTML (the common language for webpages that has many accessibility features).  Mathematics within HTML can be displayed in different ways, such as using MathJax (an engine for displaying mathematics using images that works in all web browsers) or MathML (a standard for mathematics that is supported to varying extents by different browsers).
\item PDF (Adobe portable document format with various accessibility features)
\item ePub (an electronic-book file format based on HTML)
\item Microsoft Word/Powerpoint
\item Google Docs
\end{itemize}

Since the first version of this document was produced, advances have been made in making accessible STEM PDF documents with {\LaTeX}\ \cite{latextagging}.   As such it no longer makes sense to consider HTML as the \emph{a priori} superior choice when it comes to accessibility.  That said, using HTML does not automatically mean that a document is accessible, just as using PDF for STEM documents no longer means that a document is not accessible.  Accessibility is better thought of as a spectrum, and it makes more sense to think of files as being described as ``HTML-with-accessible-features" or ``PDF-with-accessible-features".   Authors must use their tools correctly to ensure that these accessible features are available for their readers.

At the time of writing, some accessibility features within PDF are only available on certain screen readers that are only available on Windows computers.  This should be considered as a barrier to true accessibility, and hopefully the situation will improve in the near future.

Also, at present, there are technological limitations in making certain pieces of mathematics truly accessible (e.g. many kinds of commutative diagrams, or other complicated mathematical figures).  I do not try to predict what tomorrow's technology will bring, but it appears to me likely that making some of these accessible will be possible within HTML before it is possible within PDF.

\subsection{Authoring Tools}

Information about authoring tools changes quickly as technology advances.   Rather than trying to keep this current here,  I refer the reader to the American Mathematical Society resources on this topic \cite{AMS}.  Readers who are keen to make highly accessible documents with minimal author intervention or accessibility knowledge are recommended to look into PreTeXt \cite{PreTexT}.

\end{document}